\documentclass[11pt]{article}

\usepackage{bm}
\usepackage{amsfonts}
\usepackage{amssymb,amsmath,amsthm, amsfonts}

\def\sqr#1#2{{\vcenter{\vbox{\hrule height.#2pt
              \hbox{\vrule width.#2pt height#1pt \kern#1pt \vrule width.#2pt}
          \hrule height.#2pt}}}}
\def\sqr#1#2{{\vcenter{\vbox{\hrule height.#2pt
              \hbox{\vrule width.#2pt height#1pt \kern#1pt \vrule width.#2pt}
              \hrule height.#2pt}}}}
\def\3n{\negthinspace \negthinspace \negthinspace }
\def\2n{\negthinspace \negthinspace }
\def\1n{\negthinspace }

\def\={\buildrel \triangle \over =}

\def\exp{\mathop{\rm exp}}
\def\sup{\mathop{\rm sup}}
\def\inf{\mathop{\rm inf}}

\def\({\Big (}
\def\){\Big )}
\def\[{\Big[}
\def\]{\Big]}
\def\be{\begin{equation}}

\def\ee{\end{equation}}

\def\square#1{\vbox{\hrule\hbox{\vrule height#1%
     \kern#1\vrule}\hrule}}
\def\rectangle#1#2{\vbox{\hrule\hbox{\vrule height#1%
     \kern#2\vrule}\hrule}}

\font\tenbb=msbm10 \font\sevenbb=msbm7 \font\fivebb=msbm5

\newfam\bbfam
\scriptscriptfont\bbfam=\fivebb \textfont\bbfam=\tenbb
\scriptfont\bbfam=\sevenbb

\newtheorem{remark}{Remark}[section]
\newtheorem{example}{Example}[section]
\newtheorem{theorem}{Theorem}[section]

\newtheorem{proposition}{Proposition}[section]

\makeatletter

\newcommand{\Rmnum}[1]{\expandafter\@slowromancap\romannumeral #1@}
   
   \@addtoreset{equation}{section}
\makeatother

\begin{document}

\title{Comparison theorems for mean-field BSDEs whose generators depend on the law of the solution $(Y,Z)$\footnotemark[1]}
\author{ Juan Li$^{1,2,\ast\ast}$,\,\, Zhanxin Li$^{1}$,\,\, Chuanzhi Xing$^{2,\ast\ast}$ \\
   {$^1$\small  School of Mathematics and Statistics, Shandong University, Weihai,}
	{\small Weihai 264209, P. R. China.}\\
	{$^2$\small Research Center for Mathematics and Interdisciplinary Sciences,     }\\
{\small      Frontiers Science Center for Nonlinear Expectations, Ministry of Education,}\\
{\small Shandong University, Qingdao 266237, P.~R.~China.}\\
 {\small{\it E-mails: juanli@sdu.edu.cn,\,\  zhanxinli@mail.sdu.edu.cn,\,\ xingchuanzhi@sdu.edu.cn.}}
\date{June 01, 2024}
}
\renewcommand{\thefootnote}{\fnsymbol{footnote}}
\footnotetext[1]{Juan Li is supported by the NSF of China (NOs. 12031009, 11871037), the NSF of Shandong Province (NO. ZR2023ZD35), the National Key R and D Program of China (NO. 2018YFA0703900), 111 Project (No. B12023). Chuanzhi Xing is supported by the NSF of China (No. 12301177), Shandong Province (No. ZR2023QA105).\\
\noindent\ \ \  \ \ $^{\ast\ast}$Corresponding authors.}
\maketitle

\textbf{Abstract}. For general mean-field backward stochastic differential equations (BSDEs) it is well-known that we usually do not have the comparison theorem if the coefficients depend on the law of $Z$-component of the solution  process $(Y, Z)$. A natural question is whether general mean-field BSDEs whose coefficients depend on the law of $Z$ have the comparison theorem for some cases. In this paper we establish the comparison theorems for one-dimensional mean-field BSDEs whose coefficients also depend on the joint law of the solution process $(Y,Z)$. With the help of Malliavin calculus and a BMO martingale argument, we obtain two comparison theorems for different cases and a strong comparison result. In particular, in this framework, we compare not only the first component $Y$ of the solution $(Y,Z)$ for such mean-field BSDEs, but also the second component $Z$.


\section{Introduction}
Nonlinear backward stochastic differential equations (BSDEs) were first introduced by Pardoux and Peng in their pioneering paper \cite{PP1990} in 1990.
Since then the theory of BSDEs has attracted the attention of a large number of researchers because of its various applications, such as stochastic optimal control,
 stochastic differential games, mathematical finance and stochastic interpretation of partial differential equations (PDEs).

 It is worth mentioning that the comparison theorem for real-valued BSDEs turns out to be one of the classic results of BSDEs theory. The comparison theorem plays a significant role in BSDEs by providing important insights into their behavior and characteristics. It is widely used in the study of BSDEs and has profound implications for various fields such as finance, control theory, and numerical approximation.
Roughly speaking,  the comparison theorem states that if two BSDEs have the same terminal condition but with one having a larger generator than the other, then the solution to the BSDE with the larger generator will be smaller than or equal to the solution to the BSDE with the smaller generator at all time. In other words, the comparison theorem  provides a partial ordering for the solutions of BSDEs. This property allows us to compare the solutions of different BSDEs and establish relationships between them. It is particularly useful in studying the existence, uniqueness, and stability of solutions to BSDEs. In summary, the comparison theorem in BSDEs is a fundamental theorem that underlies many theoretical and practical aspects of this branch of stochastic analysis.

The comparison theorem for BSDEs were discussed the first time by Peng \cite{P1992} in order to investigate the
 viscosity solutions of second order PDEs, and then generalized by  El Karoui, Peng and Quenez \cite{EPQ1997}.
Because of crucial applications in the study of the existence of a solution of BSDEs with a continuous coefficient and
viscosity solutions of PDEs, the comparison theorem for BSDEs has received increasing attentions. For example,
the converse comparison theorem for real-valued BSDEs has been studied in \cite{CHMP2001}. In addition,
with the help of the backward stochastic viability property studied by Buckdahn, Quincampoix and Rascanu \cite{BQR2000},
Hu and Peng \cite{HP2006} studied the comparison theorem for multi-dimensional BSDEs in 2006.
It should be mentioned that in the above literature, authors only compare the first component of the solution  process $(Y,Z)$ for BSDE,
when they study the comparison theorem.

Here we are interested in BSDEs whose coefficients do not depend only on the solution processes $(Y,Z)$ but also on their joint law.
In the mathematical literature, the theory of mean field games was initially introduced by
Lasry and Lions \cite{LL2007}, simultaneously, by Huang, Caines, and Malham\'{e} \cite{HCM2006}.
Inspired by their seminal work,  Buckdahn, Djehiche, Li and Peng \cite{BDLP2009} investigated a special mean-field problem
and deduced a new type of nonlinear BSDEs which they called mean-field BSDEs. Furthermore, in 2009, Buckdahn, Li and Peng
\cite{BLP2009} obtained the existence and the uniqueness for mean-field BSDEs with Lipschitz coefficients.
When the coefficient does not depend on the law of $Z$, they gave a comparison theorem for this type of BSDE.
On the other hand, Buckdahn, Li, Peng and Rainer \cite{BLPR2017}
also considered a McKean-Vlasov SDE introduced by Kac \cite{K1956} in 1956 and the associated PDE. Subsequently, stimulated by these works, the theory of mean-field BSDEs,
as well as the associated PDEs have been studied by different authors, for example, Li \cite{L2018}, Chassagneux, Crisan and Delarue \cite{CCD2015}.
More recently, Li, Liang and Zhang \cite{LLZ2018}  studied the mean-field BSDE whose coefficient depends not only on the solution $(Y,Z)$
but also on the law $P_{Y}$ of $Y$.
In order to obtain the comparison theorem, they first introduced a weak monotonicity assumption on $f$ w.r.t. the law of $Y$ (see, (\ref{eq2.2}) in Section 2).
In fact, if $f$ is Lipschitz and monotonic, this weak monotonicity assumption  is obviously true (for details, the
readers may refer to Remark 2.1 in \cite{LLZ2018}).
With the help of the comparison result, they proved the existence theorem of such mean-field BSDE whose driving drift coefficient is only
continuous and of linear growth. Besides, under different frameworks, more and more  researchers have studied the comparison theorems of
 mean-field BSDEs, in which the coefficients depend on the solution processes but also on the law $P_{Y}$ of $Y$, such as \cite{CXZ2020,LX2022,LXP2021}.
We emphasize that the coefficients in \cite{BLP2009,CCD2015,CXZ2020,LX2022,LXP2021} are random coefficients.

In this paper, we are interested in studying the comparison theorems for the following
one-dimensional mean-field backward stochastic differential equation (BSDE):
\begin{equation}\label{eq1.1}
  \begin{split}
  Y_t=\xi+\int_t^Tf(s,Y_s,Z_s,P_{(Y_s,Z_s)})ds-\int_t^TZ_sdW_s,\ 0\leq t\leq T.
  \end{split}
\end{equation}
Let us mention that it follows from counter-examples in Buckdahn, Li and Peng \cite{BLP2009} and  Li, Liang and Zhang  \cite{LLZ2018} that
if the coefficient $f$ depends on the $Z$-component of the
law $\mu_s:=P_{(Y_s,Z_s)}$ of the solution $(Y_s,Z_s)$, we usually do not have the comparison theorem.
Therefore, for the case when the coefficient $f$ depends on the law $\mu_s$, it is difficult to get the comparison result with respect to (w.r.t.)
the first component of the solution $(Y,Z)$ for mean-field BSDE \eqref{eq1.1}, and it is even more difficult to compare  the second component
and this is even for the non mean-field case. A natural question is whether the mean-field BSDE whose generator depends on the law of $Z$ has the comparison theorem for some cases. For this, let us start with Example \ref{eg3.2}. It should be noted that in Example
\ref{eg3.2}, we get not only $Y_t^1\leq Y_t^2,$ $t\in[0,T]$, $P$-a.s., but we obtain $Z_t^1\leq Z_t^2,\ dtdP$-a.e. In other words,
 under suitable assumptions, we can compare the second component of the solution $(Y,Z)$.
On the other hand, through the calculation process, we notice that the generators in Example \ref{eg3.2} are deterministic, and the terminal conditions are Malliavin differentiable.
Therefore, this inspires us to investigate the comparison theorems for  mean-field BSDE \eqref{eq1.1} with the help of Malliavin calculus and under appropriate derivative assumptions.

In this paper, we obtain  comparison results w.r.t.
the first component of the solution $(Y,Z)$ for  mean-field BSDE \eqref{eq1.1}, but also the second component.
In fact, even already from the case of a BSDE, to our best knowledge,
there are only few works concerning the  comparison result w.r.t. the second component of the solution $(Y,Z)$ up to now.
However, when the coefficient $f$ is deterministic, we overcome this difficult, and succeed to prove the comparison theorems for
mean-field BSDE \eqref{eq1.1}.

To be more precise, our main results can be divided into two parts. In Subsection 3.1,
when the coefficient $f$ does not depend on $Z$, we give a comparison theorem for such mean-field BSDE.
Furthermore, when the coefficient $f$ depends on the solution processes $(Y,Z)$ but also on their joint law, we
also obtain a comparison theorem and a strong result in Subsection 3.2. It should be mentioned that in this case,
the terminal value $\xi$ needs to have a particular form $\xi=\Phi(W_T)$.
Let us emphasize that our results are not in conflict with the counter-examples in \cite{BLP2009} and  \cite{LLZ2018}. On the other hand, in the proof of our results, the monotonicity hypothesis
of the derivative of the terminal value (see, (H1) and (H4)) and the assumption about the derivative of the generator (see, (H3) and (H7)) are very crucial. These assumptions are new but adequate for the problem. Of course, we also give examples (see, Remark \ref{re3.1} and the Examples \ref{eg3.1} and \ref{eg3.3}) to illustrate the validity of these assumptions.

The main objective of our paper concerns the study of the comparison theorems for mean-field BSDEs under the deterministic coefficients, i.e.,
we will investigate the comparison result w.r.t. the first component of the solution $(Y,Z)$ for mean-field BSDEs, but also the second component.
Our paper is organized as follows: In Section 2 we give the necessary preparation which is needed in what follows.
In Section 3 we study the comparison theorems of solutions for one-dimensional mean-field BSDEs.
With the help of Malliavin calculus, when the coefficients do not depend on $Z$, we obtain a comparison theorem (Theorem \ref{th3.1}) in subsection 3.1.
Furthermore, when the coefficients depend on the solution processes but also on their joint law,
we also get a comparison theorem (Theorem \ref{th3.2}) and a strong result (Theorem \ref{th3.3}) in subsection 3.2.

\section{Preliminaries}
Let $T>0$ be a fixed time horizon and $(\Omega,\mathcal{F},P)$ be a given complete probability space. Let $\{W_t,\ 0\leq t\leq T\}$ be a standard Brownian motions defined on $(\Omega,\mathcal{F},P)$, with values in $\mathbb{R}^d$. We assume that there is a sub-$\sigma$-field $\mathcal{F}^0\subset\mathcal{F}$, containing all $P$-null subsets of $\mathcal{F}$, such that\\
\indent\quad(i)\ \ the  Brownian motion $W$ is independent of $\mathcal{F}^0$;

\indent\quad(ii)\ $\mathcal{F}^{0}$ is `rich enough', i.e., $\mathcal{P}_{2}(\mathbb{R}^{k})=\{P_{\vartheta},\ \vartheta\in L^{2}(\mathcal{F}^{0};\mathbb{R}^{k})\},\ k\geq1$.

\noindent Here $P_\vartheta := P\circ [\vartheta]^{-1}$ denotes the law of the random variable $\vartheta$ under the probability $P$,

\indent Let us introduce some notations and concepts, which are used frequently in what follows. Recall that $\mathcal{P}_{2}(\mathbb{R}^{k})$ is the set of the probability measures on $(\mathbb{R}^{k},\mathcal{B}(\mathbb{R}^{k}))$ with finite second moment, i.e.,  $\int_{\mathbb{R}^k}|x|^{2}\mu(dx)<\infty $. Here $\mathcal{B}(\mathbb{R}^{k})$ denotes the Borel $\sigma$-field over $\mathbb{R}^{k}$. The space $\mathcal{P}_{2}(\mathbb{R}^{k})$ be endowed with the $2$-Wasserstein metric: For $\mu,\nu\in \mathcal{P}_{2}(\mathbb{R}^k)$, we put
\begin{equation*}
W_{2}(\mu,\nu):=\inf\Big\{\Big(\int_{\mathbb{R}^k\times\mathbb{R}^k}|x-y|^{2}\gamma(dxdy)\Big)^{\frac{1}{2}}: \gamma\in\mathcal{P}_{2}(\mathbb{R}^{2k}), \gamma(.\times\mathbb{R}^k)=\mu, \gamma(\mathbb{R}^k\times.)=\nu\Big\}.
\end{equation*}

 We denote by $\mathbb{F}=(\mathcal{F}_{t})_{t\in [0,T]}$ the filtration generated by $W$
and augmented by $\mathcal{F}^0$. For $k\in \mathbb{N}$\ and\ $x,y\in\mathbb{R}^{k}$, we denote its norm and inner product, respectively,
by $|x|=\Big(\sum\limits_{i=1}^{k}x_{i}^{2}\Big)^\frac{1}{2}$, and $\langle x,y\rangle=\sum\limits_{i=1}^{k}x_{i}y_{i}.$
For $a^1,a^2\in\mathbb{R}^k$, we define $a^1\geq a^2\ \mbox{iff}\ a^1_j\geq a^2_j,\ j=1,2,\cdots,k.$
We shall also introduce the following spaces of stochastic processes: For $p\geq1$,\\
\noindent$\bullet \ L^p(\mathcal{F}_T;\mathbb{R}^d)$ is the set of $\mathcal{F}_T$-measurable random variables $ \xi: \Omega\!\rightarrow\!\mathbb{R}^d$ such that $\|\xi\|_{L^{p}}^p\!:=\!E[|\xi|^{p}]\!<\!\infty $.\\
\noindent$\bullet\ \mathcal{S}^p(0,T;\mathbb{R}^d)$ is the set of $\mathbb{F}$-adapted continuous processes  $\eta: \Omega\times[0,T]\rightarrow \mathbb{R}^d$ with $\|\eta\|_{\mathcal{S}^{p}}:=$\\
\noindent\mbox{ } \ $ \big(E\big[\sup\limits_{0\leq s\leq T}|\eta(s)|^{p}\big]\big)^{\frac{1}{p}}<\infty$.\\
\noindent$\bullet\ \mathcal{H}^p(0,T;\mathbb{R}^d)$ is the set of $\mathbb{F}$-progressively measurable processes $\eta: \Omega\times[0,T]\rightarrow \mathbb{R}^d$ with $\|\eta\|_{\mathcal{H}^{p}}:=$\\
\noindent\mbox{ } \ $\displaystyle{\big( E\big[(\int_{0}^{T}|\eta(s)|^{2}ds)^{\frac{p}{2}}\big]\big)^{\frac{1}{p}}<\infty}$.\\
\noindent$\bullet\ BMO_p$ is the set of real-valued $\mathbb{F}$-martingales $M$ such that
\begin{equation*}
||M||_{BMO_p}:=\sup_{\tau}\big\|  (E[|M_T-M_{\tau}|^p|\mathcal{F}_{\tau}])^\frac{1}{p}\big\|_{\infty}<\infty,
\end{equation*}
\noindent\mbox{ } \ where the supremum is taken over all stopping times $\tau\in[0,T]$. Here BMO martingale stands\\
\noindent\mbox{ } \ for bounded mean oscillation martingale. \\
\noindent$\bullet\ C_b^k(\mathbb{R}^p,\mathbb{R}^q)$ is the set of functions of class $C^k$ from $\mathbb{R}^p$ into $\mathbb{R}^q$ whose partial derivatives of all \\
\noindent\mbox{ } \  order less than or equal to $k$ are bounded.\\
\noindent$\bullet\ \mathcal{S}$ is the set of smooth random variables $\xi$ of the form $\xi=\varphi(W(h_1),\cdots,W(h_n))$,  $n\geq0$, with \\
\noindent\mbox{ } \ $\varphi\in C_b^{\infty}(\mathbb{R}^{n},\mathbb{R})$, $h_1,\cdots,h_n\in \mathcal{H}^2(0,T;\mathbb{R}^d)$, $W(h_i):=\displaystyle\int_0^T\langle h_i(t),dW_t\rangle$.

 Moreover, if $\xi\in S$ is of the above form, its Malliavin derivative w.r.t. $W$, denoted by $D[\cdot]$, is given by
 $ \displaystyle D_t[\xi](=( D_t^1[\xi],\cdots, D_t^d[\xi]))=\sum_{i=1}^n\frac{\partial\varphi}{\partial x_i}(W(h_1),\cdots,W(h_n))h_i(t),\ 0\leq t\leq T.$\
For $\xi\in S$, $p>1$, we define the norm $\displaystyle ||\xi||_{1,p}=\(E\[|\xi|^p+\big(\int_0^T|D_t\xi|^2dt\big)^{\frac{p}{2}}\]\)^{\frac{1}{p}}.$
From \cite{N1995} we know the operator $D$ has a closed extension to the space $\mathbb{D}^{1,p}$, the closure of $S$ with respect to the norm $||\cdot||_{1,p}$.
Observe that if $\xi\in\mathbb{D}^{1,2}$ is $\mathcal{F}_t$-measurable, $D_{\theta}\xi= 0$, $d\theta d P$-a.e., $\theta\in (t,T]$.

Recall also that, given an $\mathbb{F}$-adapted, pathwise square integrable process $\xi=(\xi_s)$,
a local martingale of the form $\displaystyle M_t=\int_{0}^{t}\xi_sdW_s$  is a BMO martingale if and
only if $\displaystyle ||M||_{BMO_2}=\sup_{\tau}\|  (E[\int_{\tau}^{T}|\xi_s|^2ds|\mathcal{F}_{\tau}])^\frac{1}{2}\|_{\infty}<\infty$.
We denote by  $\xi\cdot M$ the stochastic integral of a scalar-valued adapted
process $\xi$ with respect to a local continuous martingale $M$.
For more details of the theory of BMO martingales, the readers  may refer to Kazamaki \cite{K1994}.

\indent Let us now consider a function $f(s,y,z,\mu): [0,T]\times\mathbb{R}\times\mathbb{R}^d\times\mathcal{P}_{2}(\mathbb{R}^{1+ d})\rightarrow\mathbb{R}$
which is assumed to be $\mathbb{F}$-progressively measurable and satisfy the following standard assumptions:\\
\textbf{(A1)}  $\displaystyle\int_0^T|f(s,0,0,\delta_{0})|^2ds<+\infty$, where $\delta_{0}$ is the Dirac measure with mass at $0\in\mathbb{R}^{1+d}$.\\
\textbf{(A2)}  $f$ is Lipschitz in $(y,z,\mu)$: There exists a constant $C\!>\!0$ such that, for all $\mu,\mu'\!\in\!\mathcal{P}_{2}(\mathbb{R}^{1+ d})$,\\ \mbox{ } \ \ \ \ \ \ $y_1,y_2\in\mathbb{R}$, $s\in[0,T]$,
      $$|f(s,y_1,z_1,\mu)-f(s,y_2,z_2,\mu')|\leq C\big(W_2(\mu,\mu')+|y_1-y_2|+|z_1-z_2|\big).$$
Moreover, we assume that the terminal value $\xi$ satisfies the following standard assumption:\\
\textbf{(A3)}  $\xi\in L^{2}(\mathcal{F}_{T};\mathbb{R})$.

We consider the following general mean-field BSDEs:
\begin{equation}\label{eq2.1}
  \begin{split}
  Y_t=\xi+\int_t^Tf(s,Y_s,Z_s,P_{(Y_s,Z_s)})ds-\int_t^TZ_sdW_s,\ 0\leq t\leq T.
  \end{split}
\end{equation}
The following results on mean-field BSDEs  are by now well-known; for their proof the reader is referred to the
Theorems 2.1 and 2.2 in Li, Liang and Zhang \cite{LLZ2018}, respectively.
\begin{proposition} \label{prop2.1} (Existence and Uniqueness)
Under the assumptions (A1)-(A3), equation (\ref{eq2.1}) has a unique solution $(Y,Z)\in\mathcal{S}^2(0,T;\mathbb{R})\times\mathcal{H}^2(0,T;\mathbb{R}^d)$.
\end{proposition}

We now consider the comparison theorem for the following mean-field BSDEs when the coefficient $f$ does not depend on the law of $Z$:
\begin{equation}\label{eq2.1.1}
  \begin{split}
  Y_t=\xi+\int_t^Tf(s,Y_s,Z_s,P_{Y_s})ds-\int_t^TZ_sdW_s,\ 0\leq t\leq T.
  \end{split}
\end{equation}

\begin{proposition}(Comparison Theorem) \label{prop2.2}
Let $f_{i}:=f_{i}(s,y,z,\mu),\ i=1,2$, and the data $(\xi^{i},f_{i})$,\ $i=1,2$ satisfy (A1) and (A3). Moreover, we suppose\\
\emph{(i)} One of the $f_i$ satisfies (A2);\\
\emph{(ii)} One of the $f_i$ is such that, for all $\theta_{1},\theta_{2}\in L^{2}(\mathcal{F};\mathbb{R}),$ and all $(s,y,z)\in[0,T]\times\mathbb{R}\times\mathbb{R}^{d},$
there exists a constant $L>0$, such that
\begin{equation}\label{eq2.2}
 f(s,y,z,P_{\theta_{1}})-f(s,y,z,P_{\theta_{2}})\leq L(E[((\theta_{1}-\theta_{2})^{+})^{2}])^{\frac{1}{2}}.
\end{equation}
Let $(Y^{1},Z^{1})$ and $(Y^{2},Z^{2})$ be the solutions of mean-field BSDE (\ref{eq2.1.1}) with the data $(\xi^{1},f_{1})$ and $(\xi^{2},f_{2})$, respectively. Moreover, if $\xi^{1}\leq\xi^{2},\ P\mbox{-}a.s.$, and $f_{1}(s,y,z,\mu)\leq f_{2}(s,y,z,\mu),\ dsdP\mbox{-}a.e.,$ for all $(\mu,y,z)$, then $Y_{s}^{1}\leq Y_{s}^{2}$, for all $s\in[0,T],\ P\mbox{-}a.s.$
\end{proposition}

\section{Comparison theorems for mean-field BSDEs}
In this section, we will study  comparison theorem. For the one-dimensional comparison theorem, from the counter-examples in \cite{BLP2009} and \cite{LLZ2018} we know that, if the driver $f$ depends on the law of $Z$ or is not increasing with respect to the law of $Y$, we usually do not have the comparison theorem. On the other hand,
for a comparison theorem, we usually compare the first component of the solution $(Y,Z)$ for BSDE.

A natural question is whether the mean-field BSDE whose generator depends on the law of $Z$ can have the comparison theorem under some assumptions.
For this, let us start with the following example:
\begin{example}\label{eg3.2}
Let $d=1$ and the coefficient $f(s, y, z, P_{(\xi,\eta)})=E[\eta],$ for all $(s, y, z)\in [0, T]\times\mathbb{R}^2,\ \xi,\ \eta \in L^{2}(\mathcal{F};\mathbb{R})$. Denote the Malliavin derivative of $\xi^1$ at time $t$ by $D_t [\xi^1]$.
We consider the following mean-field BSDE:
\begin{equation*}
  Y_t^i=\xi^i+\int_t^TE[Z_s^i]ds-\int_t^TZ_s^idW_s,\ t\in[0,T],\ i=1,2.
\end{equation*}
For the equation $i=1$ we put  $\xi^1=-(W_T^+)^2$.
Then, as one can see easily,
$Z_t^1=E[D_t [\xi^1]|\mathcal{F}_t]=-2E[W_T^+|\mathcal{F}_t]$,
$E[Z_t^1]=-2\sqrt{T}E[W_1^+]=-\sqrt{\frac{2T}{\pi}},\ t\in[0,T],$
and $\displaystyle Y_t^1=E[\xi^1|\mathcal{F}_t]+\int_t^T E[Z_s^1]ds\leq0,$ $t\in[0,T]$.
For the equation $i=2$ we choose $\xi^2=0$. Then, $Y^2=0$, and $Z^2=0$. Notice that $\xi^1\leq\xi^2$, and we also have
 $Y_t^1\leq0=Y_t^2,$ $t\in[0,T]$, $P$-a.s., and $Z_t^1\leq 0=Z_t^2,\ dtdP$-a.e.
\end{example}
This example illustrates that two real-valued  mean-field BSDEs whose generators depend on the law of $Z$
may have the comparison result whenever we can compare the terminal conditions and the coefficients.
Inspired by this, we consider the following two cases to see when two real-valued  mean-field BSDEs whose generators depend on the law of $Z$
may admit the comparison result.

\subsection{The case $f(s,y,P_{(\vartheta,\eta)})$}
Now we will consider the  following mean-field BSDE:
\begin{equation}\label{eq3.1}
  \begin{split}
  Y_t=\xi+\int_t^Tf(s,Y_s,P_{(Y_s,Z_s)})ds-\int_t^TZ_sdW_s,\ 0\leq t\leq T.
  \end{split}
\end{equation}
Remark that the driving coefficient $f$ of the above BSDE depends on $Z$ only through its law.

\begin{theorem}\label{th3.1}
Let $f_{i}:=f_{i}(s,y,\mu),\ i=1,2$. We denote by $(Y^{i},Z^{i})$ the solution of mean-field BSDE (\ref{eq3.1}) with the data
$(\xi^{i},f_{i})$,\ $i=1,2$, satisfying (A1)-(A3). Moreover, we suppose\\
\textbf{\emph{(H1)}} The differentiability and a strong monotonicity of $\xi^i$:\\ \mbox{ }\ \ \ \ \ \ \
\emph{(i)} Differentiability: $\xi^i$ is Malliavin differentiable;\\ \mbox{ }\ \ \ \ \ \ \
\emph{(ii)} Monotonicity: $\xi^{1}\leq\xi^{2},\ P\mbox{-a.s.}$;\\ \mbox{ }\ \ \ \ \ \ \
\emph{(iii)}  Monotonicity of the derivative:  $(0\vee D_r^j[\xi^1])\leq D_r^j[\xi^2],\ dr d P\mbox{-a.e.},\ r\in[0,T],\ 1\leq j\leq d.$\\
\textbf{\emph{(H2)}}  $f_i$ is  monotonic in $\mu$ in the following sense:
$f_1(s,y,P_{(\vartheta,\eta)})\leq f_2(s,y,P_{(\vartheta',\eta')})$, for all $s\in[0,T]$,\\ \mbox{ } \ \ \ \ \ \
  $ y\in\mathbb{R}$, $\vartheta,\vartheta'\in L^{2}(\mathcal{F}_{T};\mathbb{R})$ with $\vartheta\leq\vartheta'$,
  $\eta,\eta'\in L^{2}(\mathcal{F}_{T};\mathbb{R}^d)$ with $\eta\leq\eta'$.\\
\textbf{\emph{(H3)}} $f_i$ is differentiable w.r.t. $y$ and the derivative $\partial_y f_i$ is bounded and monotonic:\\ \mbox{ }\ \ \ \ \ \ \
\emph{(i)} Boundedness: The derivative $\partial_y f_i$ is bounded on
$[0,T]\times\mathbb{R}\times\mathcal{P}_2(\mathbb{R}^{1\!+\!d})$;\\ \mbox{ }\ \ \ \ \ \ \
\emph{(ii)} Monotonicity in $(y,\mu)$: $(\partial_y f_1)(s,y,P_{(\vartheta,\eta)})\leq(\partial_y f_2)(s,y',P_{(\vartheta',\eta')})$,
for all $s\in[0,T]$, $ y,y'$\\ \mbox{ } \ \ \ \ \ \ \ \ \ \ \   $\in\mathbb{R}$ with $y\leq y'$, $\vartheta,\vartheta'\in$
$  L^{2}(\mathcal{F}_{T};\mathbb{R})$  with  $\vartheta\leq\vartheta'$,
$\eta,\eta'\in L^{2}(\mathcal{F}_{T};\mathbb{R}^d)$  with  $\eta\leq\eta'$.

\noindent Then $Y_t^1\leq Y_t^2,\ t\in[0,T],$\ $P$-a.s., and $Z_t^1\leq Z_t^2,\ dtdP$-a.e.
\end{theorem}
\begin{proof}
For simplicity of redaction but without loss of generality, let us restrict to dimension $d = 1$.
Under the assumptions (A1)-(A3), it follows from Proposition \ref{prop2.1} that the mean-field BSDE (\ref{eq3.1})
with the data $(\xi^{1},f_{1})$ and $(\xi^{2},f_{2})$ have a unique solution $(Y^{1},Z^{1})$ and $(Y^{2},Z^{2})$ in $\mathcal{S}^2(0,T;\mathbb{R})\times\mathcal{H}^2(0,T;\mathbb{R})$, respectively.
Let $(Y^{i,0},Z^{i,0}):=(0,0)$, $i=1,2$.  We consider the following BSDE:
\begin{equation}\label{eq3.2}
  Y_t^{i,n}=\xi^i+\int_t^T f_i(s,Y_s^{i,n},P_{(Y_s^{i,n-1},Z_s^{i,n-1})})ds-\int_t^T Z_s^{i,n} dW_s,\ t\in[0,T],\ n\geq1,\ i=1,2.
\end{equation}
Again from Proposition \ref{prop2.1}, we know  BSDE (\ref{eq3.2}) has a unique solution $(Y^{i,n},Z^{i,n}),\ n\geq1,$ for $i=1,2$.
Furthermore, from standard arguments, it can easily be verified that $(Y^{i,n},Z^{i,n})_{n\geq1}$ is a Cauchy sequence in $\mathcal{S}^2(0,T;\mathbb{R})\times\mathcal{H}^2(0,T;\mathbb{R})$, which means
\begin{equation*}
  E[\sup\limits_{t\in[0,T]}|Y_t^i-Y_t^{i,n}|^2]+E[\int_0^T|Z_s^i-Z_s^{i,n}|^2 ds]\rightarrow0,\ n\rightarrow\infty.
\end{equation*}
We now suppose that, for some $n\geq1,$\\
\textbf{(A.n-1)}\quad $ Y_t^{1,n-1}\leq Y_t^{2,n-1},\ t\in[0,T],\ P\mbox{-a.s.},\quad Z_t^{1,n-1}\leq Z_t^{2,n-1},\ dtdP$-a.e.

Observe that this holds, in particular, for $n=1$. Let us show that (A.n-1) implies (A.n), for $n\geq1$.
From (A1)-(A3), (H1)-(ii) and (H2), we deduce with the help of the comparison theorem-Proposition \ref{prop2.2}  that
\begin{equation}\label{eq3.3}
Y_t^{1,n}\leq Y_t^{2,n},\ t\in[0,T],\ P\mbox{-a.s.}
\end{equation}
Define $\theta_s^{i,n}:=(s,Y_s^{i,n},P_{(Y_s^{i,n-1},Z_s^{i,n-1})})$. We take now the Malliavin derivative in \eqref{eq3.2}: For $r\leq t$,
\begin{equation}\label{eq3.4}
D_r[Y_t^{i,n}]=D_r[\xi^i]+\int_t^T(\partial_y f_i)(\theta_s^{i,n})D_r[Y_s^{i,n}]ds-\int_t^T D_r[Z_s^{i,n}]dW_s,\ t\in[0,T].
\end{equation}
From (H1)-(iii) we know $D_r[\xi^2]\geq0,\ dr d P\mbox{-a.e.},\ r\in[0,T]$. Therefore, it follows from a by now classical argument that
\begin{equation}\label{eq3.5}
D_r[Y_t^{2,n}]\geq0,\  dr d P\mbox{-a.e.},\ 0\leq r\leq t\leq T.
\end{equation}
Moreover, for $\overline{Y}^n:=Y^{2,n}-Y^{1,n},\ \overline{Z}^n:=Z^{2,n}-Z^{1,n},$ and $\overline\xi:=\xi^2-\xi^1,$
\begin{equation}\label{eq3.6}
\begin{split}
D_r[\overline{Y}_t^n]=&\ D_r[\overline\xi]+\int_t^T\big((\partial_y f_2)(\theta_s^{2,n})-(\partial_y f_1)(\theta_s^{1,n})\big)D_r[Y_s^{2,n}]ds+\int_t^T(\partial_y f_1)(\theta_s^{1,n})D_r[\overline{Y}_s^n]ds\\
&\ -\int_t^T D_r[\overline{Z}_s^n]dW_s,\ 0\leq r\leq t\leq T.
\end{split}
\end{equation}
Making use of (H1)-(iii), (H3), (A.n-1), (\ref{eq3.3}) and (\ref{eq3.5}), we obtain
$$\big((\partial_y f_2)(\theta_s^{2,n})-(\partial_y f_1)(\theta_s^{1,n})\big)D_r[Y_s^{2,n}]\geq0.$$
Thus, as also $D_r[\overline\xi]\geq0$, from a standard BSDE estimate we obtain $D_r[Y_t^{2,n}]-D_r[Y_t^{1,n}]=D_r[\overline{Y}_t^n]\geq0,\ 0\leq r\leq t\leq T.$
Letting $r\leq t\downarrow r$, this yields $Z_t^{2,n}-Z_t^{1,n}=\overline{Z}_t^n\geq0,\ t\in[0,T].$ This proves (A.n).
Passing to the limit, then the proof is complete.
\end{proof}

\begin{remark}\label{re3.1}
Observe that unlike Theorem 2.2 in \cite{LLZ2018}, here we  get  the comparison theorem with respect to the solution $(Y,Z)$, when the driver $f$ depends on the solution processes $Y$ but also on their joint law of $(Y,Z)$. It should be mentioned that Theorem \ref{th3.1} is not in contradiction with the counter-example (Example 2.1) in \cite{LLZ2018}.
In particular, letting the dimension $d=1$, we discuss a particular case, i.e., we consider the following mean-field BSDE:
\begin{equation*}
  Y_t^i=\xi^i+\int_t^T f_i(P_{Z_s^i})ds-\int_t^T Z_s^i dW_s,\ t\in[0,T],\ i=1,2.
\end{equation*}
Here we see, in particular that
$D_r[Y_t^i]=D_r[\xi^i]-\int_t^TD_r Z_s^i dW_s,\ r\leq t\leq T$, and
$ Z_r^i=D_r[\xi^i]-\int_r^T D_r[Z_s^i]dW_s, $ i.e., $Z_r^i=E[D_r[\xi^i]|\mathcal{F}_r],\ drdP$-a.e.
Obviously, $D_r[\xi^1]\leq D_r[\xi^2],\ r\in[0,T],$ implies here $Z_r^1\leq Z_r^2,\ drdP$-a.e.,
and so
\begin{equation*}
Y_t^1 =E[\xi^1+\int_t^T f_1(P_{Z_s^1})ds|\mathcal{F}_t]\leq E[\xi^2+\int_t^T f_2(P_{Z_s^2})ds|\mathcal{F}_t]=Y_t^2,\ t\in[0,T].
\end{equation*}
Here we have used that $f_1(P_\eta)\leq f_2(P_{\eta'}),$ for $\eta,\eta'\in L^{2}(\mathcal{F}_{T};\mathbb{R})$ with $ \eta\leq\eta'$,
 following from (H2). This property is essential. We also give a counter-example here.
\end{remark}

\begin{example}\label{eg3.1}
The coefficient $f_i(s,y,P_{(\vartheta,\eta)})=f_i(P_\eta)=E[|\eta|],\ i=1,2,$ does, obviously, not satisfy (H2)$;$
and we don't have the comparison theorem result for it.
Indeed, let $d=1$. We consider the following mean-field BSDE:
\begin{equation*}
  Y_t^i=\xi^i+\int_t^TE[|Z_s^i|]ds-\int_t^TZ_s^idW_s,\ t\in[0,T],\ i=1,2.
\end{equation*}
Let us put $\xi^1=-(W_T^+)^2$. Then, obviously,
$Z_t^1=E[D_t [\xi^1]|\mathcal{F}_t]=-2E[W_T^+|\mathcal{F}_t]$,
$E[|Z_t^1|]=2\sqrt{T}E[W_1^+]=\sqrt{\frac{2T}{\pi}},\ t\in[0,T],$
and $\displaystyle Y_0^1=E[\xi^1]+\int_0^T E[|Z_s^1|]ds=-\frac{T}{2}+\sqrt{\frac{2}{\pi}}T^{\frac{3}{2}}>0,$ for $T\geq1$.
We choose $\xi^2=0$, and so we have $Y^2=0$, $Z^2=0$, in particular, $Y_0^2=0$.
Therefore, for $T\geq1$,\ $Y_0^2=0<Y_0^1,$ although $\xi^1\leq0=\xi^2,$ and $D_r[\xi^1]\leq0=D_r[\xi^2],\ r\in[0,T].$
\end{example}

\begin{remark}\label{re3.2}
The conclusion from (A.n-1) to (A.n) can also be made by using the multi-dimensional BSDE comparison result by Hu and Peng \cite{HP2006} for the BSDE
$$
\Big(\!\begin{array}{c}
Y_t^{i,n}\\
D_r[Y_t^{i,n}]\\
\end{array}\!\Big)\!=\!\Big(\!\begin{array}{c}
\xi^i\\
D_r[\xi^i]\\
\end{array}\!\Big)\!+\!\int_t^T\!\!\Big(\!\begin{array}{c}
f_i(s,Y_s^{i,n},P_{(Y_s^{i,n-1},Z_s^{i,n-1})})\\
\partial_y f_i(s,Y_s^{i,n},P_{(Y_s^{i,n-1},Z_s^{i,n-1})})D_r[Y_s^{i,n}]\\
\end{array}\!\Big)ds\!-\!\int_t^T\!\!\Big(\!\begin{array}{c}
Z_s^{i,n}\\
D_r[Z_s^{i,n}]\\
\end{array}\!\Big)dW_s,
$$
with $0\leq r\leq t\leq T,\ i=1,2$.
\end{remark}

\subsection{The case $f(s,y,z,P_{(\vartheta,\eta)})$}
Let us now suppose for the terminal value $\xi$ that it is of the form $\xi=\Phi(W_T)$, where $\Phi(x):\mathbb{R}\rightarrow\mathbb{R}$.
That is, we consider the following mean-field BSDE:
\begin{equation}\label{eq3.7}
  \begin{split}
  Y_t=\Phi(W_T)+\int_t^Tf(s,Y_s,Z_s,P_{(Y_s,Z_s)})ds-\int_t^TZ_sdW_s,\ 0\leq t\leq T.
  \end{split}
\end{equation}
For the study of this BSDE we restrict ourselves to dimension  $d = 1$.

\begin{theorem}\label{th3.2}
Let $f_{i}:=f_{i}(s,y,z,\mu),\ i=1,2$. We denote by $(Y^{i},Z^{i})$ the solution of mean-field BSDE (\ref{eq3.7}) with the data
$(\Phi_i(W_T),f_{i})$,\ $i=1,2$ satisfying (A1)-(A3). Moreover, we suppose\\
\textbf{\emph{(H4)}} The differentiability and strong monotonicity of $\Phi_i$:\\ \mbox{ }\ \ \ \ \ \ \
\emph{(i)} Differentiability:  $\Phi_1$ is differentiable w.r.t. $x$, $\Phi_2$ is twice differentiable w.r.t. $x$;\\ \mbox{ }\ \ \ \ \ \ \
\emph{(ii)} Monotonicity: $\Phi_1(x)\leq\Phi_2(x)$, for all $ x\in\mathbb{R}$;\\ \mbox{ }\ \ \ \ \ \ \
\emph{(iii)}  Monotonicity of the derivative: $0\leq\partial_x\Phi_1(x)\leq\partial_x\Phi_2(x)$, for all $ x\in\mathbb{R}$.\\
\textbf{\emph{(H5)}} The boundedness of the derivatives of $\Phi_2$:
 There exists $C>0$ such that, for all $ x\in\mathbb{R}$,
  $$|\partial_x\Phi_2(x)|\leq C\ \mbox{and}\ 0\leq\partial_x^2\Phi_2(x)\leq C.$$
\textbf{\emph{(H6)}}  $f_i$ is  monotonic in $\mu$: $f_1(s,y,z,P_{(\vartheta,\eta)})\leq f_2(s,y,z,P_{(\vartheta',\eta')})$, for all $s\in[0,T]$, $ y, z\in\mathbb{R}$,\\ \mbox{ } \ \ \ \ \ \
$\vartheta,\vartheta',\eta,\eta'\in L^{2}(\mathcal{F}_{T};\mathbb{R})$
with $\vartheta\leq\vartheta'$, $\eta\leq\eta'$.\\
\textbf{\emph{(H7)}} $f_1$ is differentiable w.r.t. $(y,z)$, $f_2$ is twice differentiable w.r.t. $(y,z)$, and  these derivatives\\ \mbox{ } \ \ \ \ \ \   have the following properties: \\ \mbox{ } \ \ \ \ \
\emph{(i)} Boundedness: All derivatives of
 order 1 and 2 are uniformly bounded over $[0,T]\times\mathbb{R}\times$\\ \mbox{ } \ \ \ \ \ \ \ \ \ \ $\mathbb{R}\times\mathcal{P}_{2}(\mathbb{R}^{2})$;\\ \mbox{ }\ \ \ \ \ \
\emph{(ii)} Convexity: $(y,z)\rightarrow f_2(y,z,P_{(\vartheta,\eta)})$ is convex;\\ \mbox{ }\ \ \ \ \ \
\emph{(iii)} Lipschitz continuity:  $z\rightarrow\partial_y f_2(s,y,z,P_{(\vartheta,\eta)})z$ is Lipschitz, uniformly in $(s,y,P_{(\vartheta,\eta)})$;\\ \mbox{ }\ \ \ \ \ \
\emph{(iv)} Monotonicity in $(y,\mu)$:  For all
$ y,y'\!\in\!\mathbb{R}$ with $y\leq y'$,
 $\vartheta,\vartheta',\eta,\eta'\!\in\!  L^{2}(\mathcal{F}_{T};\mathbb{R})$ with $\vartheta\leq\vartheta'$, \\ \mbox{ } \ \ \ \ \ \ \ \ \ \ \
$\eta\leq\eta'$,  $s\in[0,T]$,  $z\in\mathbb{R}$,
 \begin{align*}
(a)\ &(\partial_y f_1)(s,y,z,P_{(\vartheta,\eta)})\!\leq\!(\partial_y f_2)(s,y',z,P_{(\vartheta',\eta')});\\
(b)\ &(\partial_z f_1)(s,y,z,P_{(\vartheta,\eta)})\!\leq\!(\partial_z f_2)(s,y',z,P_{(\vartheta',\eta')}).
 \end{align*}
Then $Y_t^1\leq Y_t^2,\ t\in[0,T],$\ $P$-a.s., and $Z_t^1\leq Z_t^2,\ dtdP$-a.e.
\end{theorem}

\begin{proof}
Under the assumptions (A1)-(A3), it follows from Proposition \ref{prop2.1} that the mean-field BSDE (\ref{eq3.7})
with the data $(\Phi_1(W_T),f_{1})$ and $(\Phi_2(W_T),f_{2})$ has a unique solution $(Y^{1},Z^{1})$ and $(Y^{2},Z^{2})$ in $\mathcal{S}^2(0,T;\mathbb{R})\times\mathcal{H}^2(0,T;\mathbb{R})$, respectively.
Let $(Y^{i,0},Z^{i,0}):=(0,0)$, $i=1,2$.  We consider the following BSDE:
\begin{equation}\label{eq3.8}
  Y_t^{i,n}\!=\!\Phi_i(W_T)\!+\!\int_t^T\!\! f_i(s,Y_s^{i,n},Z_s^{i,n},P_{(Y_s^{i,n-1},Z_s^{i,n-1})})ds\!-\!\int_t^T\!\! Z_s^{i,n} dW_s,\ t\!\in\![0,T],\ n\geq1,\ i=1,2.
\end{equation}
Thanks to Proposition \ref{prop2.1} BSDE (\ref{eq3.8}) has a unique solution $(Y^{i,n},Z^{i,n}),\ n\geq1,\ i=1,2$.
Furthermore, by standard arguments, it can easily be verified that $(Y^{i,n},Z^{i,n})_{n\geq1}$ is a Cauchy sequence in $\mathcal{S}^2(0,T;\mathbb{R})\times\mathcal{H}^2(0,T;\mathbb{R})$,
\begin{equation*}
  E[\sup\limits_{t\in[0,T]}|Y_t^i-Y_t^{i,n}|^2]+E[\int_0^T|Z_s^i-Z_s^{i,n}|^2 ds]\rightarrow0,\ n\rightarrow\infty.
\end{equation*}
We now suppose that, for some $n\geq1,$\\
\textbf{(A.n-1)}\quad $ Y_t^{1,n-1}\leq Y_t^{2,n-1},\ t\in[0,T],\ P\mbox{-a.s.},\quad Z_t^{1,n-1}\leq Z_t^{2,n-1},\ dtdP$-a.e.

Notice that this holds, in particular, for $n=1$. Let us show that (A.n-1) implies (A.n), for $n\geq1$.
From (A1)-(A3), (H4)-(ii), (H6) and the comparison theorem-Proposition \ref{prop2.2}, we obtain
\begin{equation}\label{eq3.9}
Y_t^{1,n}\leq Y_t^{2,n},\ t\in[0,T],\ P\mbox{-a.s.}
\end{equation}
Define $\theta_s^{i,n}:=(s,Y_s^{i,n},Z_s^{i,n},P_{(Y_s^{i,n-1},Z_s^{i,n-1})})$. We take now the Malliavin derivative in \eqref{eq3.8};
standard arguments show its existence, and that for $0\leq r\leq t\leq T$, $i=1,2$,
\begin{equation}\label{eq3.10}
D_r[Y_t^{i,n}]\!=\!\partial_x\Phi_i(W_T)\!+\!\int_t^T\!\!\big((\partial_y f_i)(\theta_s^{i,n})D_r[Y_s^{i,n}]
\!+\!(\partial_z f_i)(\theta_s^{i,n})D_r[Z_s^{i,n}]\big)ds\!-\!\int_t^T\!\! D_r[Z_s^{i,n}]dW_s.
\end{equation}
From the uniqueness of the solution it follows that $D_r[Y_t^{i,n}]=D_{r'}[Y_t^{i,n}],\ r\vee r'\leq t\leq T,$
and $D_r[Z_t^{i,n}]=D_{r'}[Z_t^{i,n}],\ dtdP$-a.e. on $[r\vee r',T]\times\Omega.$
Thus, $Z_t^{i,n}=P$-$\lim\limits_{r\uparrow t}D_r[Y_t^{i,n}]=D_r[Y_t^{i,n}],$ for all $r\leq t.$
As $\partial_x\Phi_i(x)\geq0$, $x\in\mathbb{R}$ (see (H4)-(iii)), a standard argument involving Girsanov transformation shows
\begin{equation}\label{eq3.10+1}
Z_t^{i,n}=D_r[Y_t^{i,n}]\geq0,\   0\leq r\leq t\leq T,\ i=1,2.
\end{equation}
We  put $\widetilde{Z}_s^{i,n}:= D_r[Z_s^{i,n}],\ 0\leq r\leq s\leq T.$
Then, from \eqref{eq3.10},
\begin{equation}\label{eq3.11}
Z_t^{i,n}=\partial_x\Phi_i(W_T)+\int_t^T\big((\partial_y f_i)(\theta_s^{i,n})Z_s^{i,n}
   +(\partial_z f_i)(\theta_s^{i,n})\widetilde{Z}_s^{i,n}\big)ds-\int_t^T\widetilde{Z}_s^{i,n}dW_s,\ t\in[0,T].
\end{equation}
Let us set $f_i^n(s,\cdot,\cdot):= f_i(s,\cdot,\cdot,P_{(Y_s^{i,n-1},Z_s^{i,n-1})}).$ Then the system \eqref{eq3.8}, \eqref{eq3.11} writes:
\begin{equation}\label{eq3.12}
\begin{split}
Y_t^{i,n}=\Phi_i(W_T)+\int_t^T f_i^n(s,Y_s^{i,n},Z_s^{i,n})ds-\int_t^T Z_s^{i,n}dW_s,\ t\in[0,T],\\
\end{split}
\end{equation}
and
\begin{equation}\label{eq3.13}
\begin{split}
Z_t^{i,n}=&\ \partial_x\Phi_i(W_T)\!+\!\int_t^T\!\!\!\big((\partial_y f_i^n)(s,Y_s^{i,n},D_r[Y_s^{i,n}])Z_s^{i,n}+(\partial_z f_i^n)(s,Y_s^{i,n},D_r[Y_s^{i,n}])\widetilde{Z}_s^{i,n}\big)ds\\
&\ -\int_t^T\widetilde{Z}_s^{i,n}dW_s,\ 0\leq r\leq t\leq T,\ i=1,2.\\
\end{split}
\end{equation}
Obviously, the solution $(Z^{i,n},\widetilde{Z}^{i,n})$ of BSDE \eqref{eq3.13} has the following property:
\begin{equation}\label{eq3.13+1}
  |Z_t^{i,n}|\leq C, \  \mbox{for all}\ t\in[0,T],\ P\mbox{-a.s.},\ \mbox{and}\ \widetilde{Z}^{i,n}\in\mathcal{H}^p(0,T;\mathbb{R}),\ p\geq 1,\ i=1,2.
\end{equation}
Thus,  $Z^{i,n}\cdot W(=D_r[Y^{i,n}]\cdot W)$ is a BMO martingale. On the other hand, we notice that \eqref{eq3.11} can also be rewritten as follows:
\begin{equation}\label{eq3.13-new1}
\begin{split}
Z_t^{i,n}=&\ \partial_x\Phi_i(W_T)\!+\!\int_t^T\!\!\!\big((\partial_y f_i^n)(s,Y_s^{i,n},D_r[Y_s^{i,n}])D_r[Y_s^{i,n}]+(\partial_z f_i^n)(s,Y_s^{i,n},D_r[Y_s^{i,n}])\widetilde{Z}_s^{i,n}\big)ds\\
&\ -\int_t^T\widetilde{Z}_s^{i,n}dW_s,\ 0\leq r\leq t\leq T,\ i=1,2.\\
\end{split}
\end{equation}
From Lemma 2.1 in \cite{HT2016}, we have $\widetilde{Z}^{i,n}\cdot W$ is a BMO martingale.
Differentiating \eqref{eq3.11} formally in the sense of Malliavin we get for $i=1, 2$,
\begin{equation}\label{eq3.14}
\begin{split}
D_r&[Z_t^{i,n}]= \partial_x^2\Phi_i(W_T)+\int_t^T\Big((\partial_y^2 f_i)(\theta_s^{i,n})(Z_s^{i,n})^2+(\partial_y\partial_z f_i)(\theta_s^{i,n})(Z_s^{i,n}D_r[Z_s^{i,n}]+Z_s^{i,n}\widetilde{Z}_s^{i,n})\\
&\ +(\partial_z^2 f_i)(\theta_s^{i,n})\widetilde{Z}_s^{i,n}D_r[Z_s^{i,n}]\Big)ds+\int_t^T\Big((\partial_y f_i)(\theta_s^{i,n})D_r[Z_s^{i,n}]+(\partial_z f_i)(\theta_s^{i,n})D_r[\widetilde{Z}_s^{i,n}]\Big)ds\\
&\ -\int_t^T D_r[\widetilde{Z}_s^{i,n}]dW_s,\ 0\leq r\leq t\leq T.\\
\end{split}
\end{equation}
We remark that this BSDE is linear, for $0\leq r\leq t\leq T$,
\begin{equation}\nonumber
\begin{split}
D_r[Z_t^{i,n}]= \partial_x^2\Phi_i(W_T)+\int_t^T\Big(h_{1,i}(s)+h_{2,i}(s)D_r[Z_s^{i,n}]+h_{3,i}(s)D_r[\widetilde{Z}_s^{i,n}]\Big)ds
 -\int_t^T D_r[\widetilde{Z}_s^{i,n}]dW_s,
\end{split}
\end{equation}
where the coefficients
\begin{equation}\nonumber
\begin{split}
h_{1,i}(s)& :=(\partial_y^2 f_i)(\theta_s^{i,n})(Z_s^{i,n})^2+(\partial_y\partial_z f_i)(\theta_s^{i,n})Z_s^{i,n}\widetilde{Z}_s^{i,n},\\
h_{2,i}(s)& :=(\partial_y\partial_z f_i)(\theta_s^{i,n})Z_s^{i,n}+(\partial_z^2 f_i)(\theta_s^{i,n})\widetilde{Z}_s^{i,n}+(\partial_y f_i)(\theta_s^{i,n}),\\
h_{3,i}(s)& :=(\partial_z f_i)(\theta_s^{i,n}),\ s\in[0,T].
\end{split}
\end{equation}
Notice that $h_{1,i}\in\mathcal{H}^p(0,T;\mathbb{R})$, $h_{2,i}\in\mathcal{H}^p(0,T;\mathbb{R})$, $p\geq2$,
while $h_{3,i}$ is bounded and $\mathbb{F}$-adapted, $i=1,2$.
Recall that $\widetilde{Z}^{i,n}\cdot W$ is a BMO martingale. Thus from Theorem 10 in \cite{BC2008}, \eqref{eq3.14} has a unique solution $(D_r[Z_t^{i,n}],D_r[\widetilde{Z}_s^{i,n}])_{0\leq r\leq t\leq T}\in\mathcal{S}^p(0,T;\mathbb{R})\times\mathcal{H}^p(0,T;\mathbb{R})$, $p\geq2$.
The existence and the uniqueness of this solution allows to use standard arguments to show that
it is really the Malliavin derivative.
From the uniqueness of the solution $(D_r[Z_s^{2,n}],D_r[\widetilde{Z}_s^{2,n}])$ we have
$$D_r[Z_t^{2,n}]=D_{r'}[Z_t^{2,n}],\ r\vee r'\leq t\leq T,\ P\mbox{-a.s.};\quad
D_r[\widetilde{Z}_t^{2,n}]=D_{r'}[\widetilde{Z}_t^{2,n}],\ dtdP\mbox{-a.e. on }[r\vee r',T]\times\Omega.$$
Consequently, $\widetilde{Z}_t^{2,n}=P$-$\lim\limits_{t>r\uparrow t}D_r[Z_t^{2,n}]=D_r[Z_t^{2,n}],\ 0\leq r\leq t\leq T.$ And we also introduce the notation
$$ \widehat{Z}_t^{2,n}:= D_r[\widetilde{Z}_t^{2,n}],\ (0\leq r\leq)t\leq T. $$
Thus, (\ref{eq3.14}) takes the form
\begin{equation}\label{eq3.15}
\begin{split}
\widetilde{Z}_t^{2,n}=&\ \partial_x^2\Phi_2(W_T)+\int_t^T\Big\langle\Big(\begin{array}{cc}
\partial_y^2 f_2 & \partial_y\partial_z f_2\\
\partial_y\partial_z f_2 & \partial_z^2 f_2\\
\end{array}\Big)(\theta_s^{2,n})\Big(\begin{array}{c}
Z_s^{2,n}\\
\widetilde{Z}_s^{2,n}\\
\end{array}\Big),\Big(\begin{array}{c}
Z_s^{2,n}\\
\widetilde{Z}_s^{2,n}\\
\end{array}\Big)\Big\rangle ds\\
&\ +\int_t^T\big((\partial_y f_2)(\theta_s^{2,n})\widetilde{Z}_s^{2,n}
+(\partial_z f_2)(\theta_s^{2,n})\widehat{Z}_s^{2,n}\big)ds-\int_t^T\widehat{Z}_s^{2,n}dW_s,\ t\in[0,T],\ i=1,2.
\end{split}
\end{equation}
Denote $\zeta_s^{2,n}:=\Big\langle\Big(\begin{array}{cc}
\partial_y^2 f_2 & \partial_y\partial_z f_2\\
\partial_y\partial_z f_2 & \partial_z^2 f_2\\
\end{array}\Big)(\theta_s^{2,n})\Big(\begin{array}{c}
Z_s^{2,n}\\
\widetilde{Z}_s^{2,n}\\
\end{array}\Big),\Big(\begin{array}{c}
Z_s^{2,n}\\
\widetilde{Z}_s^{2,n}\\
\end{array}\Big)\Big\rangle$, then we deduce from (H7)-(i), (ii) that
\begin{equation}\label{eq3.16}
  0\leq\zeta_s^{2,n}\leq C\big(|Z_s^{2,n}|^2+|\widetilde{Z}_s^{2,n}|^2\big),\ \mbox{where}\ C> 0\ \mbox{ may vary line to line}.
\end{equation}
In order to use Girsanov transformation, we define the probability measure $\widetilde{P}_n$ and the process $\widetilde{W}^{n}$ as follows
\begin{equation*}
  \begin{split}
    &\widetilde{P}_n\!:= \!\mbox{exp}\Big\{\int_0^T\!\! (\partial_z f_2)(\theta_s^{2,n})dW_s
         \!-\!\frac{1}{2}\int_0^T\!\!|(\partial_z f_2)(\theta_s^{2,n})|^2ds\Big\}\!\cdot\! P,\quad\!
     \widetilde{W}_t^n\!:=\!W_t\!-\!\int _0^t\!\!(\partial_z f_2)(\theta_s^{2,n})ds,\ t\!\in\![0,T].
  \end{split}
\end{equation*}
Due to Girsanov's theorem  $\widetilde{W}^n=(\widetilde{W}_t^{n})_{0\leq t\leq T}$ is an $(\mathbb{F},\widetilde{P}_n)$-Brownian motion.\\

Let $\displaystyle \alpha_{t,r}^n:=\exp\Big\{$ $\displaystyle \int_t^r(\partial_y f_2)(\theta_s^{2,n})ds\Big\}$,\ $0\leq t\leq r\leq T$. Then we get
\begin{equation*}
 \widetilde{Z}_t^{2,n}=\alpha_{t,T}^n\cdot\partial_x^2\Phi_2(W_T)
+\int_t^T\alpha_{t,s}^n\cdot\zeta_s^{2,n}ds-\int_t^T\alpha_{t,s}^n\cdot\widehat{Z}_s^{2,n}d\widetilde{W}_s^n,\ t\in[0,T].
\end{equation*}
Hence, $\displaystyle\widetilde{Z}_t^{2,n}=\widetilde{E}^n\big[\alpha_{t,T}^n\cdot\partial_x^2\Phi_2(W_T)
+\int_t^T\alpha_{t,s}^n\cdot\zeta_s^{2,n}ds|\mathcal{F}_t\big]$,
where $\displaystyle \widetilde{E}^n[\cdot]=\int_{\Omega}(\cdot)d\widetilde{P}_n$.
Furthermore, it follows from (H5), (H7)-(i), \eqref{eq3.13+1}, \eqref{eq3.16} and the fact that $\widetilde{Z}^{2,n}\cdot W$ is a BMO martingale that
\begin{equation}\label{eq3.16+1}
\begin{split}
   0\leq&\ \widetilde{Z}_t^{2,n}\leq C\Big(1+\widetilde{E}^n\big[\int_t^T(|Z_s^{2,n}|^2+|\widetilde{Z}_s^{2,n}|^2)ds|\mathcal{F}_t\big]\Big)\leq C,\ t\in[0,T],
\end{split}
\end{equation}
 which allows to go back to (\ref{eq3.13}).
Let $ \overline{\Phi}:=\Phi_2-\Phi_1,\ \overline{Z}^n:= Z^{2,n}-Z^{1,n},\ \widetilde{Z}^n:=\widetilde{Z}^{2,n}-\widetilde{Z}^{1,n}. $
Then we have the following equation:
\begin{equation}\label{eq3.17}
\begin{split}
\overline{Z}_t^n=&\ \partial_x\overline{\Phi}(W_T)\!+\!\int_t^T\!\!\Big(
\big((\partial_y f_2^n)(s,Y_s^{2,n},Z_s^{2,n})Z_s^{2,n}-(\partial_y f_1^n)(s,Y_s^{1,n},Z_s^{1,n})Z_s^{1,n}\big)\\
&\ +(\partial_z f_1^n)(s,Y_s^{1,n},Z_s^{1,n})\widetilde{Z}_s^n+\big((\partial_z f_2^n)(s,Y_s^{2,n},Z_s^{2,n})\!-\!(\partial_z f_1^n)(s,Y_s^{1,n},Z_s^{1,n})\big)\widetilde{Z}_s^{2,n}
\Big)ds\\
&\ -\int_t^T\widetilde{Z}_s^n dW_s,\ t\in[0,T].
\end{split}
\end{equation}
As due to (A.n-1),  \eqref{eq3.9}, \eqref{eq3.10+1} and
(H7)-(iii), (iv)-(a) imply:
\begin{equation}\label{eq3.18}
(\partial_y f_2^n)(s,Y_s^{2,n},Z_s^{2,n})Z_s^{2,n}-(\partial_y f_1^n)(s,Y_s^{1,n},Z_s^{1,n})Z_s^{1,n}\geq\gamma_s^n\overline{Z}_s^n,\ s\in[0,T],
\end{equation}
for some bounded progressively measurable process $\gamma^n$.
We note that, since $\partial_z^2 f_2$ is bounded, $z\rightarrow(\partial_z f_2)(s,y,z,P_{(\vartheta,\eta)})$ is Lipschitz,
uniformly w.r.t. $(s,y,P_{(\vartheta,\eta)})$. Thus, from (H7)-(iv)-(b), (A.n-1) and \eqref{eq3.16+1},
we see that, for some bounded, progressively measurable process $\widetilde{\gamma}^n,$
\begin{equation}\label{eq3.19}
\begin{split}
&\ \big((\partial_z f_2^n)(s,Y_s^{2,n},Z_s^{2,n})-(\partial_z f_1^n)(s,Y_s^{1,n},Z_s^{1,n})\big)\widetilde{Z}_s^{2,n}\\
=&\ \big((\partial_z f_2^n)(s,Y_s^{2,n},Z_s^{2,n})-(\partial_z f_2^n)(s,Y_s^{2,n},Z_s^{1,n})\big)\widetilde{Z}_s^{2,n}\\
&\ +\big((\partial_z f_2^n)(s,Y_s^{2,n},Z_s^{1,n})-(\partial_z f_1^n)(s,Y_s^{1,n},Z_s^{1,n})\big)\widetilde{Z}_s^{2,n}\\
\geq&\ -C|(\partial_z f_2^n)(s,Y_s^{2,n},Z_s^{2,n})-(\partial_z f_2^n)(s,Y_s^{2,n},Z_s^{1,n})|=\widetilde{\gamma}_s^n\overline{Z}_s^n,\ s\in[0,T].
\end{split}
\end{equation}
The estimates (\ref{eq3.18}) and (\ref{eq3.19}) make it standard to conclude from (\ref{eq3.17}) that
$$ Z_t^{2,n}-Z_t^{1,n}=\overline{Z}_t^n\geq0,\ t\in[0,T]. $$
Therefore, $Y^{2,n}\geq Y^{1,n}$ and $Z^{2,n}\geq Z^{1,n},$ i.e., (A.n) holds true.
By induction we have (A.n) for all $n\geq1.$ Consequently, taking the limit, as $n\rightarrow\infty,$ we obtain
$Y_t^1\leq Y_t^2,\ t\in[0,T],$\ $P$-a.s., and $Z_t^1\leq Z_t^2,\ dtdP$-a.e.
\end{proof}

\begin{example}\label{eg3.3}
Let $f_i(s,y,z,P_{(\vartheta,\eta)})=-y+z+E[\vartheta]+E[\eta],\ i=1,2$. It can easily be checked that $f_i,\ i=1,2$
satisfy the assumptions (A1), (A2), (H6) and (H7). We consider the following mean-field BSDE:
\begin{equation}\label{eq3.23+1}
  Y_t^i=\xi^i+\int_t^T\Big(-Y_s^i+Z_s^i+E[Y_s^i]+E[Z_s^i]\Big)ds-\int_t^TZ_s^idW_s,\ t\in[0,T],\ i=1,2.
\end{equation}
Let $\varepsilon>0$ and let us consider the function $U_\varepsilon\in C(\mathbb{R})$:
\begin{equation*}
{U_\varepsilon(y):=}
\left\{ \begin{aligned}
        \ 0,\qquad\qquad y&\leq0;\\
        \ \frac{1}{\varepsilon^2}y,\qquad\quad  0&\leq y\leq \varepsilon;\\
        \ \frac{2}{\varepsilon}-\frac{1}{\varepsilon^2}y,\quad \varepsilon&\leq y\leq 2\varepsilon;\\
        \ 0,\qquad\qquad  y&\geq2\varepsilon.
\end{aligned} \right.
\end{equation*}
We define
\begin{equation*}
{R_\varepsilon(y):=\int^y_{-\infty} U_\varepsilon(z)dz=}
\left\{ \begin{aligned}
        \ 0,\qquad\qquad\qquad\qquad y&\leq0;\\
        \ \frac{1}{2\varepsilon^2}y^2,\qquad\qquad\qquad\ 0&\leq y\leq \varepsilon;\\
        \ \frac{2}{\varepsilon}y-\frac{1}{2\varepsilon^2}y^2-1,\quad\quad \varepsilon&\leq y\leq 2\varepsilon;\\
        \ 1,\qquad\qquad\qquad\qquad y&\geq2\varepsilon;
\end{aligned} \right.
\end{equation*}
and $\displaystyle \Phi_2(x):=\int^{x}_{-\infty} R_\varepsilon(y)dy$, for any $x\in \mathbb{R}$.\\

We take $\xi^1=0$, then, obviously, the unique solution \eqref{eq3.23+1} is $(Y^1,Z^1)=(0,0)$.
On the other hand, we take $\xi^2=\Phi_2(W_T)$. Then, from \eqref{eq3.23+1}, for $r\leq t$,
\begin{equation}\nonumber
  D_r[Y_t^2]=R_\varepsilon(W_T)+\int_t^T(-D_r[Y_s^2])ds-\int_t^TD_r[Z_s^2](dW_s-ds),\ t\in[0,T],
\end{equation}
i.e.,
\begin{equation}\nonumber
Z_t^2=R_\varepsilon(W_T)+\int_t^T(-Z_s^2)ds-\int_t^TD_r[Z_s^2](dW_s-ds),\ t\in[0,T].
\end{equation}
Consequently,
$Z_t^2=E\big[\exp\big\{-\frac{3}{2}(T-t)+W_T-W_t\big\}\cdot R_\varepsilon(W_T)|\mathcal{F}_t\big]\geq0$, $t\in[0,T]$,
and thus, $E[Z_t^2]=E\big[\exp\big\{-\frac{3}{2}(T-t)+W_T-W_t\big\}\cdot R_\varepsilon(W_T)]\geq0$, $t\in[0,T]$. Moreover, also from \eqref{eq3.23+1}
$\displaystyle E[Y_t^2]=E\big[\xi^2+\int_t^T2E[Z_s^2]ds\big]\geq0$, $t\in[0,T]$.
Finally, again from  \eqref{eq3.23+1}, with a standard argument for linear BSDEs,
 $\displaystyle Y_t^2=E\Big[\exp\big\{-\frac{3}{2}(T-t)+W_T-W_t\big\}\cdot\xi^2+\int_t^T \exp\big\{-\frac{3}{2}(s-t)+W_s-W_t\big\}\cdot\big(E[Y_s^2]+E[Z_s^2]\big)ds\big|\mathcal{F}_t\Big]\geq0,$ $t\in[0,T]$.
Notice that the assumptions in Theorem 3.2 are now all satisfied.
Moreover, we observe that we have $Y_t^1=0\leq Y_t^2,$ $t\in[0,T]$, $P$-a.s., and $Z_t^1=0\leq Z_t^2,\ dtdP$-a.e.
\end{example}

Finally, we give still a strong comparison result.

\begin{theorem}\label{th3.3}(A strong comparison result)
Let the data $(\Phi_1(W_T),f_{1})$ and $(\Phi_2(W_T),f_{2})$ satisfy the assumptions in Theorem \ref{th3.2}.
 Denote by $(Y^{1},Z^{1})$ and $(Y^{2},Z^{2})$ the solutions of mean-field BSDE (\ref{eq3.7}) with the data $(\Phi_1(W_T),f_{1})$
 and $(\Phi_2(W_T),f_{2})$, respectively. Then, for all $0\leq t_0<T$,  $P$-a.s. on $\{Y_{t_0}^2=Y_{t_0}^1\},$
 \begin{equation*}
   \begin{split}
     \emph{(i)}&\ Y_s^2=Y_s^1,\ s\in[t_0,T], \mbox{and } Z_s^2=Z_s^1,\ ds\mbox{-a.e. on } [t_0,T];\\
     \emph{(ii)}&\ f_2(Y_s^2,Z_s^2,P_{(Y_s^2,Z_s^2)})=f_1(Y_s^1,Z_s^1,P_{(Y_s^1,Z_s^1)}),\ ds\mbox{-a.e. on }[t_0,T].
    \end{split}
 \end{equation*}
\end{theorem}
\begin{proof}
From \eqref{eq3.7}, with the notation $\overline{\xi}=\xi^2-\xi^1,\ \overline{Y}=Y^2-Y^1,\ \overline{Z}=Z^2-Z^1,$ and observing that there are two bounded, progressively measurable processes $\alpha$ and $\beta$ such that
$$ f_2(Y_s^2,Z_s^2,P_{(Y_s^2,Z_s^2)})-f_2(Y_s^1,Z_s^1,P_{(Y_s^2,Z_s^2)})=\alpha_s\overline{Y}_{\!s}+\beta_s\overline{Z}_s,\ s\in[0,T], $$
we get, for $t\in[0,T]$,
\begin{equation}\label{eq3.20}
\begin{split}
\overline{Y}_{\!t}\!=\!\overline{\xi}\!+\!\int_t^T\!\!\!\Big(\!\big(f_2(Y_s^1,Z_s^1,P_{(Y_s^2,Z_s^2)})\!-\!f_1(Y_s^1,Z_s^1,P_{(Y_s^1,Z_s^1)})\big)
\!+\!\alpha_s\overline{Y}_s\!+\!\beta_s\overline{Z}_s\!\Big)ds\!-\!\int_t^T\!\overline{Z}_sdW_s.
\end{split}
\end{equation}
Consequently,  for $\displaystyle\widetilde{P}=L_T P,$ with $L_T=\exp\{\int_0^T\beta_s dW_s-\frac{1}{2}\int_0^T\beta_s^2 ds\},$
\begin{equation*}
  \begin{split}
   \overline{Y}_{\!t_0}=\widetilde{E}\Big[e^{\int_{t_0}^T\alpha_s ds}\overline{\xi}+\int_{t_0}^T e^{\int_{t_0}^s\alpha_r dr}(f_2(Y_s^1,Z_s^1,P_{(Y_s^2,Z_s^2)})-f_1(Y_s^1,Z_s^1,P_{(Y_s^1,Z_s^1)}))ds|\mathcal{F}_{t_0}\Big].
  \end{split}
\end{equation*}
Thus, on $\{\overline{Y}_{\!\!t_0}=0\},\ \overline{\xi}=0$ and $f_2(Y_s^1,Z_s^1,P_{(Y_s^2,Z_s^2)})=f_1(Y_s^1,Z_s^1,P_{(Y_s^1,Z_s^1)}),\ dsdP$-a.e. on $[t_0,T]\times\{\overline{Y}_{\!\!t_0}=0\}.$ Then, from \eqref{eq3.20},
$$ \overline{Y}_{\!t}=0,\ t\in[t_0,T],\ \overline{Z}_t=0,\ dt\mbox{-a.e. on }[t_0,T],\ P\mbox{-a.s. on }\{\overline{Y}_{\!\!t_0}=0\}. $$
The statement follows.
\end{proof}

\end{document}